\documentclass{amsproc}

\usepackage{cleveref}
\usepackage{amssymb}
\usepackage{tikz-cd}
\usepackage{enumitem}

\newcommand{\SL}{\textup{SL}}
\newcommand{\Q}{\mathbb{Q}}
\newcommand\C{\mathbb{C}}
\newcommand\Z{\mathbb{Z}}
\newcommand{\Qbar}{{\overline{\Q}}}
\newcommand{\mm}{\mathfrak{m}}
\newcommand{\F}{\mathbb{F}}

\DeclareMathOperator{\Aut}{Aut} 
\DeclareMathOperator{\Pic}{Pic}
\DeclareMathOperator{\Gal}{Gal}
\DeclareMathOperator{\End}{End}

\newtheorem{theorem}{Theorem}[section]
\newtheorem{lemma}[theorem]{Lemma}

\theoremstyle{definition}
\newtheorem{definition}[theorem]{Definition}

\newtheorem*{ack}{Acknowledgements}
\newtheorem{question}[theorem]{Open Question}

\theoremstyle{remark}
\newtheorem{remark}[theorem]{Remark}

\numberwithin{equation}{section}

\begin{document}

\title{Modularity over $\C$ implies modularity over $\Q$}

\author{\sc Barinder S. Banwait}
\address{Barinder S. Banwait \\
Lodha Mathematical Sciences Institute, Mumbai, India}
\email{barinder.s.banwait@gmail.com}
\urladdr{https://barinderbanwait.github.io/}


\begin{abstract}
We give an account of Mazur's proof that, for an elliptic curve over $\Q$, if it admits a nonconstant mapping from $X(N)$ defined over the complex numbers $\C$, for some $N$, then it also admits a nonconstant mapping from $X_0(M)$ for some (possibly different) $M$ defined over the rational numbers $\Q$. We also briefly discuss two open questions of Khare concerning uniformisations of elliptic curves by noncongruence modular curves.

This expository note is based on the author's expository talk in March 2022 during the 2\textsuperscript{nd} Trimester Program on \emph{Modularity and the Generalised Fermat Equation} held online and organised by Bhaskaracharya Pratishthana in Pune, India.
\end{abstract}

\maketitle

\section{Introduction}

The multiple different versions of the ``Modularity theorem'' that one may state can be daunting for the beginner to navigate, and in this particular journey the book by Diamond and Shurman \cite{diamond2005first} offers a detailed and lucid map of terrain. To quote their preface: \emph{...this book's aim is not to prove Modularity but to state its different versions, showing some of the relations among them and how they connect to different areas of mathematics.}

The question arises then of where to start this journey? That is, what is the absolute simplest formulation of the Modularity theorem? 

The first version of Modularity that appears in the main text of \cite{diamond2005first} is the following complex analytic version.

\begin{theorem}[Modularity Theorem, Version $X_0(N)_\C$]
    Let $E$ be a complex elliptic curve with $j(E) \in \Q$. Then for some positive integer $N$ there exists a surjective holomorphic function of compact Riemann surfaces:
    \[ X_0(N) \to E. \]
\end{theorem}

Mazur, however, in his wonderful expository article \cite{mazur1991number} aimed at explaining the statement of (what was then known as the Taniyama-Shimura-Weil conjecture but which of course now is) the Modularity theorem to a general audience, takes a slightly different approach. While still keeping the definition geometric, his approach is to not single out any particular congruence subgroup, and to allow $E$ to admit a dominant holomorphic function from \emph{any} modular curve associated to any congruence subgroup, what he calls a \textbf{hyperbolic uniformisation of arithmetic type}. In effect, this yields the following version. 

\begin{theorem}[Modularity Theorem, Version $X(N)_\C$]
    Let $E$ be a complex elliptic curve with $j(E) \in \Q$. Then for some positive integer $N$ there exists a surjective holomorphic function of compact Riemann surfaces:
    \[ X(N) \to E. \]
\end{theorem}

However, as noted by both Diamond and Shurman as well as Mazur, the version of Modularity that has direct arithmetic applications, including to the resolution of Fermat's Last Theorem (an observation made by Frey \cite{frey1986links}), is the following version in the language of arithmetic geometry.

\begin{theorem}[Modularity Theorem, Version $X_0(N)_\Q$]
    Let $E$ be an elliptic curve over $\Q$. Then for some positive integer $N$ there exists a surjective morphism over $\Q$ of algebraic curves over $\Q$:
    \[ X_0(N) \to E. \]
\end{theorem}

By base-change, this latter version $X_0(N)_\Q$ implies the former versions $X_0(N)_\C$ and $X(N)_\C$. The main goal of this expository note is to show that all three are in fact equivalent. That $X(N)_\C$ implies $X_0(N)_\Q$ is stated without proof by Diamond and Shurman, with the reader being referred to the `Technical Appendix' of Mazur's expository article for the details. We will give a summary of this proof in \Cref{sec:compare_mazur}, after explaining the proof that version $X_0(N)_\C$ implies $X_0(N)_\Q$ in \Cref{sec:proof_of_converse}, which may be considered one step in Mazur's proof (the Eichler-Shimura step). This proof is essentially contained in the proof of Theorem 8.8.2 of \cite{diamond2005first}. We conclude in \Cref{sec:khare} with some open questions of Khare concerning \emph{Belyi uniformisations} of elliptic curves by more general finite-index subgroups of $\SL_2(\Z)$, not necessarily congruence subgroups.

\begin{remark}
The above version $X_0(N)_\Q$ of modularity merely asserts the existence of such an $N$ for which one has a map $X_0(N) \to E$. There is an even stronger version of modularity where $N$ is taken to be the conductor of $E$. That these two versions of modularity are equivalent is a very deep result of Carayol \cite[Theorem 1]{carayol1986representations}. See the discussion after the proof of Theorem 8.8.2 in \cite{diamond2005first} for a sketch of some of the ideas involved in this.
\end{remark}

\ack{
I thank Devendra Tiwari for the invitation to give a talk in the 2\textsuperscript{nd} Trimester Program on \emph{Modularity and the Generalised Fermat Equation} held online and organised by Bhaskaracharya Pratishthana in Pune, India, from which this note arose. I am very grateful to Alex Bartel for patiently and at length explaining to me Serre's tensor construction and twisting abelian varieties by Artin representations the evening before my talk, as well as for hosting me at Glasgow University while giving the online talk; to Neil Dummigan for spotting a mistake in a previous version of the manuscript; to Barry Mazur for comments on an earlier version of the manuscript; to Bjorn Poonen for carefully explaining the material presented in \Cref{sec:C_to_Qbar}, and for insisting on greater details in \Cref{sec:compare_mazur}; to Dipendra Prasad, the chair of the talk, for asking pertinent questions related to quotients of $J_1(N)$ versus $J_0(N)$; and to Armand Brumer, for helping answer them. I also thank Peter Gr\"{a}f, Robert Pollack, and David Rohrlich for helpful discussions during the preparation of this note, as well as the anonymous referee for a careful reading of the previous version, for spotting errors in some of the proofs, and for making many helpful suggestions.

This work is supported by grant 550023 from the Simons Foundation for the collaboration `Arithmetic Geometry, Number Theory, and Computation'.
}

\section{Specialising from $\C$ to $\overline{\Q}$}\label{sec:C_to_Qbar}

In both \Cref{sec:proof_of_converse,sec:compare_mazur}, we will descend a map of curves, a priori defined over $\Qbar$, from $\C$ to $\Qbar$. We isolate this step into the following result.

\begin{lemma}\label{lem:descent_C_to_Qbar}
Let $X$ and $Y$ be proper curves defined over $\overline{\Q}$. If there exists a nonconstant morphism $X_\C \to Y_\C$, then there exists a nonconstant morphism $X \to Y$ defined over $\overline{\Q}$.
\end{lemma}

This is a standard result that we learned from Bjorn Poonen. The usual proof of this goes via specialisation as follows; it is a variant of what is discussed in \cite[Section 3.2.2]{poonen2017rational}.

\begin{proof}
    Suppose that we have a nonconstant morphism $f \colon X_\C \to Y_\C$.
Let $x \in X(\Qbar)$.
There are at most $\deg f$ elements of the infinite set $X(\Qbar)$ that have the same image in $Y$ as $x$ does.
Choose $x' \in X(\Qbar)$ such that $f(x') \ne f(x)$.
There is a finitely generated $\Qbar$-algebra $R$ such that $X_\C \to Y_\C$ arises as the base extension of a morphism $X_R \to Y_R$ of $R$-schemes, whose generic fiber is a nonconstant morphism.
Adjoin to $R$ the inverses of the differences of the coordinates of $f(x')$ and $f(x)$.
For any maximal ideal $\mm$ of $R$, the quotient $R/\mm$ is $\Qbar$, by the weak Nullstellensatz.
Reducing $X_R \to Y_R$ modulo $\mm$ gives a morphism $X \to Y$,
and it maps $x$ and $x'$ to distinct points, since $f(x')$ and $f(x)$ are distinct mod $\mm$.
\end{proof}

The following is an alternative argument.

\begin{proof}
If the genus of $Y$ is greater than $1$, then the set of nonconstant morphisms $X_\C \to Y_\C$ is finite, and they are thus all defined over $\Qbar$, since otherwise applying elements of $\Aut(\C/\Qbar)$ would produce infinitely many morphisms.

Now suppose that $g(Y)=1$, so $Y$ is an elliptic curve.
The set of morphisms $X \to Y$ is the same as the set of rational maps $X \dashrightarrow Y$ (by the valuative criterion for properness), which is the same as $Y(k(X))$.
The function field Mordell--Weil theorem (known as the Lang--N\'eron theorem) says that $Y(k(X))/Y(k)$ is finitely generated.
Equivalently, if we fix $x \in X$, then the set of morphisms $X \to Y$ mapping $x$ to $0$ is a finitely generated abelian group.
In particular, it is countable.
All of these morphisms must be defined over $\Qbar$, because again otherwise we would obtain uncountably many morphisms by applying elements of $\Aut(\C/\Qbar)$.
\end{proof}

\section{Version $X_0(N)_\C$ implies version $X_0(N)_\Q$: Diamond and Shurman's proof}\label{sec:proof_of_converse}

As stated in the introduction, this proof follows the proof of Theorem 8.8.2 of \cite{diamond2005first}.

\begin{proof}
Let $E/\Q$ be an elliptic curve. Write $N_E$ for its conductor. By assumption, there exists some $N \geq 1$ (not necessarily $N_E$) and a surjective holomorphic function of compact Riemann surfaces:
\begin{equation}
\alpha_\C : X_0(N) \to E.
\end{equation}

By \Cref{lem:descent_C_to_Qbar} we obtain a map $\alpha_\Qbar : X_0(N) \to E$ of algebraic curves over $\Qbar$, which induces a map between abelian varieties:
\[ \alpha_{\Qbar,\ast} : J_0(N)_{/\Qbar} \to E_{/\Qbar}.\]
However, we have that $J_0(N)$ is $\Q$-isogenous to a direct sum of $\Q$-simple abelian varieties $A_f$ associated to (Galois equivalence classes of) newforms,
\[ J_0(N) \to \bigoplus_f (A_{f})^{m_f}. \]
Here the sum is taken over a set of representatives $f \in S_2(\Gamma_0(M_f))$ at levels $M_f$ dividing $N$, and each $m_f$ is the number of divisors of $N/M_f$. (See \cite[Theorem 6.6.6]{diamond2005first} for this result in the context of $J_1(N)$, and the ensuing discussion there, as well as in the proof of Theorem 8.8.2, for the analogous decomposition for $J_0(N)$. The multiplicities $m_f$ arise essentially from the observation that if $f$ is a newform of level $M | N$, then the set of $f(n\tau)$, for $n$ running through the divisors of $N/M$, are linearly independent.) In what follows we will use rather the dual isogeny (and extension of scalars to $\Qbar$):

\[ \bigoplus_f (A_{f/\Qbar})^{m_f} \to J_0(N)_{/\Qbar}. \]

We now consider the following diagram, for $p \nmid N_EN$:

\[
\begin{tikzcd}
\bigoplus_f (A_{f/\Qbar})^{m_f} \arrow{r}{\prod_f\left(a_p(f) - a_p(E)\right)^{m_f}} \arrow{d} &[6em] \bigoplus_f (A_{f/\Qbar})^{m_f} \arrow{d}\\
J_0(N)_{/\Qbar} \arrow{r}{T_p - a_p(E)} & J_0(N)_{/\Qbar} \arrow{r}{\alpha_{\Qbar,\ast}} & E_{/\Qbar}.
\end{tikzcd}
\]

To define the top map we have used the fact that $A_f$ admits an action of $a_p(f)$ acting as the Hecke operator $T_p$ at the level of $A_f$; see \cite[Diagram 6.15]{diamond2005first} and the surrounding discussion for more details.

We claim that this diagram has the following properties:

\begin{enumerate}[label=(\alph*)]
    \item
    If, for some $f$, there exists a $p$ such that $a_p(f) \neq a_p(E)$, then the top map restricted to $(A_{f/\Qbar})^{m_f}$ is surjective.
    \item
    The square commutes.
    \item
    The composite map on the bottom is zero.
\end{enumerate}

Given these properties, a simple diagram chase allows us to conclude: if for each $f$ there exists a $p \nmid N_EN$ such that $a_p(f) \neq a_p(E)$, then chasing around the diagram from top left to bottom right, we map to zero, contradicting the surjectivity of $\alpha_{\Qbar,\ast}$. Thus, there must exist a newform $f \in S_2(\Gamma_0(M_f))$ such that $a_p(f) = a_p(E)$ for almost all $p$. This means that the $\ell$-adic Galois representations of $A_f$ and $E$ are isomorphic, and hence, by Faltings's Isogeny Theorem, are $\Q$-isogenous; the desired map then arises as follows:
\[ X_0(N) \hookrightarrow J_0(N) \to A_f \to E. \]
It remains to establish the above properties (a), (b) and (c). For (a), recall that $a_p(f)$ is an algebraic integer, and $a_p(E)$ is a rational integer. Thus their difference $\delta$ satisfies a minimal monic polynomial with rational integer coefficients,
\[ \delta^e + a_1\delta^{e-1} + \cdots + a_{e-1}\delta + a_e = 0 \]
with $a_e \neq 0$. The resulting relation $\delta(\delta^{e-1} + a_1\delta^{e-2} + \cdots + a_{e-1}) = -a_e$ shows that $\delta$ is surjective on $(A_{f/\Qbar})^{m_f}$.

Property (b) is a standard result arising from the decomposition of $J_0(N)$ into the direct sum of the $A_f$.

For property (c), we consider the following diagram, with the bottom row above now becoming the top row:

\begin{equation}\label{diag:main}
\begin{tikzcd}
J_0(N)_{/\overline{\Q}} \arrow{r}{T_p - a_p(E)} \arrow{d} &[3em] J_0(N)_{/\overline{\Q}} \arrow{r}{\alpha_{\ast}} \arrow{d} &[3em] E_{/\overline{\Q}} \arrow{d}\\
\widetilde{J_0}(N)_{/\overline{\F_p}} \arrow{r}{\sigma_{p,\ast} + \sigma_p^\ast - a_p(E)} \arrow{d}{1} &[3em] \widetilde{J_0}(N)_{/\overline{\F_p}} \arrow{r}{\tilde{\alpha}_{\ast}} &[3em] \widetilde{E}_{/\overline{\F_p}} \arrow{d}{1}\\
\widetilde{J_0}(N)_{/\overline{\F_p}} \arrow{r}{\tilde{\alpha}_{\ast}} &[3em] \widetilde{E}_{/\overline{\F_p}} \arrow{r}{\sigma_{p,\ast} + \sigma_p^\ast - a_p(E)}  &[3em] \widetilde{E}_{/\overline{\F_p}}
\end{tikzcd}
\end{equation}

Here the vertical arrows from the top row to the middle row indicate specialisations to the special fibre in characteristic $p$, and $\sigma_{p,\ast}, \sigma_p^\ast$ are respectively the natural induced pushforward and pullback maps from $\widetilde{J_0}(N)_{/\overline{\F_p}} \cong \Pic^0(\widetilde{X_0(N)})$ to itself, induced from the Frobenius automorphism $\sigma_p : \overline{\F_p} \to \overline{\F_p}$ taking $x \mapsto x^p$ (see \cite[\S 8.2]{diamond2005first} for more details of these maps).

As a morphism from a variety to a curve, the top row is either surjective or is constant; so we need only show that the top row is not surjective. We omit the details for why this diagram commutes (they can be consulted at the end of the proof of Theorem 8.8.2 in \cite{diamond2005first}). The second map on the bottom row is zero by \cite[Proposition 8.3.2]{diamond2005first}; thus the bottom row is zero, and hence so too is the middle row. Thus the top row followed by the vertical right arrow (reduction mod $p$) is zero; but the reduction map is surjective, meaning that the top row can't be surjective, as desired.
\end{proof}

\section{Version $X(N)_\C$ implies version $X_0(N)_\Q$: Mazur's proof}\label{sec:compare_mazur}

\begin{proof}
Let $E/\Q$ be an elliptic curve. By assumption, there exists $N \geq 1$ and a surjective holomorphic function of compact Riemann surfaces:
\begin{equation}\label{eq:original_morphism}
X(N) \to E.
\end{equation}

We proceed in a series of small steps.

\subsection*{Step 1. Reduce to the non-CM case.} If $E$ has complex multiplication over $\overline{\Q}$, then this was proved by Shimura \cite{shimura1971elliptic}. The proof (involving Gr\"{o}\ss encharaktere) uses ideas and techniques somewhat different from the rest of this note so we shall omit it here.

Henceforth, we assume that $E$ does not have CM over $\overline{\Q}$.

\subsection*{Step 2. Specialise map~\ref{eq:original_morphism} from $\C$ to $\overline{\Q}$ to a finite Galois extension $L/\Q$.} Specialising from $\C$ to $\overline{\Q}$ follows from \Cref{lem:descent_C_to_Qbar}. Specialising further to $L$ is clear because there can be only finitely many polynomials that define the curve.

\subsection*{Step 3. Replace $X(N)$ with $X_1(N^2)$.} This arises from considering the \textbf{conjugate congruence subgroup} ($\Gamma = \Gamma(N)$):
\[ \tilde{\Gamma} := \left( \begin{array}{cc}
    N & 0 \\
    0 & 1
\end{array} \right)^{-1} \Gamma \left( \begin{array}{cc}
    N & 0 \\
    0 & 1
\end{array} \right) = \Gamma_0(N^2) \cap \Gamma_1(N) \supseteq \Gamma_1(N^2).
\]
We thus get a map
\[ X_1(N^2) \to X(\tilde{\Gamma}); \]
but the point is that there is an isomorphism
\begin{align*}
    X(\tilde{\Gamma}) &\longrightarrow X(\Gamma)\\
    z &\longmapsto Nz,
\end{align*}
so we get an $L$-rational map
\[ X_1(N)_{/L} \to E_{/L} \]
(after possibly changing $N$ to $N^2$) as desired, for some finite Galois extension $L/\Q$.

\subsection*{Step 4. Pass to the Jacobian $J_1(N)$ and decompose.} The map from the previous step induces an $L$-rational map between abelian varieties
\[ J_1(N)_{/L} \to E_{/L}. \]
However, we have \cite[Theorem 6.6.6]{diamond2005first} that $J_1(N)$ is $\Q$-isogenous to a direct sum of $\Q$-simple abelian varieties $A_f$ associated to (Galois equivalence classes of) newforms,
\[ J_1(N) \to \bigoplus_f A_f^{m_f}. \]
Here the sum is taken over a set of representatives $f \in S_2(\Gamma_1(M_f))$ at levels $M_f$ dividing $N$, and each $m_f$ is the number of divisors of $N/M_f$. 

We thus obtain the map
\[ \bigoplus_f (A_{f,/L})^{m_f} \to E_{/L}. \]

\subsection*{Step 5. Take the Weil Restriction of scalars down to $\Q$.} By Serre's tensor construction \cite[\S III.1.3]{serre1995cohomologie}, we obtain a map defined over $\Q$
\begin{equation}\label{eq:main}
\bigoplus_f A_{f}^{m_f} \otimes_\Z \Q[G] \to E \otimes_\Z \Q[G],
\end{equation}
where $G := \Gal(L/\Q)$, and $\Q[G]$ is the regular representation. See \cite{milne1972arithmetic} or \cite[Section 4]{mazur2007twisting} for the details of the relation between the Weil restriction of scalars, and twists by Artin representations.

\subsection*{Step 6. Identify, up to $\Q$-isogeny, some $A_f$ with a twist of $E$ by an irreducible representation $M$ of $G$.} On the left-hand side of \Cref{eq:main}, each $A_f$ arises as a simple factor, because $\Q[G]$ contains the trivial representation. 

If it were the case that every $A_f$ were killed by the map in \Cref{eq:main}, then the composite map $J_1(N) \to J_1(N) \otimes_\Z \Q[G] \to E \otimes_\Z \Q[G]$ would be zero, which we can see is not the case as follows. We base change to $L$, where we have $J_1(N)_{/L} \to (J_1(N)_{/L})^d \to (E_{/L})^d$, where $d=[L:\Q]$. The first map here is the diagonal embedding, so by surjectivity of $J_1(N)_{/L} \to E_{/L}$ the composite cannot be 0. 

There is therefore some $A_f$ whose image under the map is nonzero. Restricting the map to this $A_f$ yields an isogeny onto its image, which is an isogeny factor of the right-hand side of \Cref{eq:main}. This right-hand side breaks up (up to isogeny) as a product of twists of $E$ by the irreducible representations of $G$. \emph{A priori} these twists need not be simple over $\Q$; however here they are, since $E$ does not have CM \cite[Proposition 2.6]{bartel2013simplicity}. We therefore obtain, for this $A_f$, an irreducible representation $M_f$ of $G$ such that
\[ A_f \sim_\Q E \otimes_\Q M_f. \]
Morally, at this point, one wants to ``untwist $E$'' to obtain
\[ A_f \otimes_\Q M_f^{-1} \sim_\Q E, \]
and then to relate the left-hand side here with the abelian variety associated to a newform (ultimately to be the newform $\psi$ in Step 9). Doing this rigorously is essentially the content of the rest of the proof. 

In the sequel, for ease of notation, we replace $M_f$ by $M$.

\subsection*{Step 7. Show that $M$ is an $F$-vector space of dimension $1$, where $F$ is the Fourier coefficient field of $f$.} Since the dimension of $A_f$ is $[F : \Q]$, and the dimension of $E \otimes_\Q M$ is $\dim_\Q(M)$, we obtain 
\[ [F : \Q] = \dim_\Q(M).\]
We also have from \cite[Corollary 4.2]{ribet1980twists} that
\[ F = \End_\Q^0(A_f).\] 
This gives us an embedding 
\[ F = \End_\Q^0(A_f) = \End_\Q^0(E \otimes_\Q M) \hookrightarrow \End_\Q^0((E_{/L}) \otimes_\Q M) = \End_\Q(M), \]
with the last equality coming from \cite[Lemma 2.1]{bartel2013simplicity}. This establishes $M$ as an $F$-vector space of dimension $1$.

\subsection*{Step 8. Show that $G$ acts on $M$ via a character $\chi$.} The previous step established that $(E_{/L}) \otimes_\Q M$ admits an $F$-action. But by definition of the Serre tensor construction, it also admits a $G$-action. These actions commute, since the $F$-action on $M$ is $\Q$-rational. Therefore $G$ acts $F$-linearly on $M$. Since $\dim_F(M) = 1$, we obtain that the $G$-action may be expressed via a character
\[ \chi : G \to F^\times. \]

\subsection*{Step 9. Identify $\chi$ as a Dirichlet character.} Since $F^\times$ is abelian, $\chi$ factors through the abelianization $G^{ab}$ of $G$. By the Kronecker-Weber theorem we obtain a surjection $\Gal(\Q(\zeta_n)/\Q) \twoheadrightarrow G^{ab}$ for some $n$, and hence a Dirichlet character \[(\Z/n\Z)^\times \cong \Gal(\Q(\zeta_n)/\Q) \twoheadrightarrow G^{ab} \to F^\times \subset \C^\times.\] By a standard abuse of notation, in the sequel we will refer to this map as $\chi$.

\subsection*{Step 10. Show $a_p(f) = a_p(E)\chi(p)$ for almost all $p$.} We apply the Eichler-Shimura relations to $A_f$; i.e. we consider the diagram analogous to \ref{diag:main} (as Diamond and Shurman do later in Section 8.8, and removing the several notational decorations for greater clarity):
\[
\begin{tikzcd}
A_f \arrow{r}{a_p(f) - a_p(E)\chi(p)} \arrow{d} &[6em] A_f \arrow{r}{\alpha} \arrow{d} &[6em] E \otimes M \arrow{d}\\
\widetilde{A_f} \arrow{r}{\sigma_{p,\ast} + \widetilde{\langle p\rangle}_{\ast}\sigma_p^\ast - a_p(E)\chi(p)} \arrow{d}{1} &[6em] \widetilde{A_f} \arrow{r}{\tilde{\alpha}} &[6em] \widetilde{E \otimes M} \arrow{d}{1}\\
\widetilde{A_f} \arrow{r}{\tilde{\alpha}} &[6em] \widetilde{E \otimes M} \arrow{r}{\sigma_{p,\ast} + \widetilde{\langle p\rangle}_{\ast}\sigma_p^\ast - a_p(E)\chi(p)}  &[6em] \widetilde{E \otimes M}.
\end{tikzcd}
\]
Here $\alpha$ is the isogeny from Step 6, $p$ is a prime of good reduction for both $A_f$ and $E \otimes M$ (i.e. $p$ does not divide the level of $f$, the conductor of $E$, or the modulus of $\chi$), the top vertical maps are reduction mod $p$, and $\langle p\rangle$ is the $p$\textsuperscript{th} Diamond operator.

If $a_p(f) \neq a_p(E)\chi(p)$ then the top row is surjective, as is $E \otimes M \to \widetilde{E \otimes M}$; but the bottom row is zero, and hence the middle row is zero, giving a contradiction because all the rectangles commute (the top-left one commuting is precisely the relation of Eichler-Shimura). We thus obtain, for almost all primes $p$, that
\[ a_p(f) = a_p(E)\chi(p). \]

\subsection*{Step 11. Obtain a newform whose Fourier coefficients agree with $a_p(E)$ for almost all $p$.} The relation from the previous step implies that, for the twisted modular form $f \otimes \bar{\chi}$, we have $a_p(f \otimes \bar{\chi}) = a_p(E)$ for almost all $p$. By a standard corollary of Strong Multiplicity One (see \cite[Proposition 5.8.4]{diamond2005first}) we obtain a newform $\psi$ with Fourier coefficient field $\Q$ such that $a_p(\psi) = a_p(E)$ for almost all $p$.

\subsection*{Step 12. Conclude via Galois representations and the Isogeny theorem.} We thus have that the $\ell$-adic Galois representations attached to $A_\psi$ and $E$ are isomorphic for all $\ell$. By Faltings's Isogeny Theorem (which in this context of modular abelian varieties was proven earlier by Ribet \cite[Theorem 6.1]{ribet1980twists}), this means that we have a $\Q$-isogeny
\[ A_\psi \sim_\Q E. \]  

Moreover, the isomorphism of $\ell$-adic Galois representations means in particular that the characteristic polynomials of Frobenius-at-$p$ must agree for almost all $p$, which shows (by considering the constant term of these quadratic polynomials) that the Nebentypus of $\psi$ is trivial. We therefore deduce a dominant morphism over $\Q$ as follows:
\[ X_0(N) \hookrightarrow J_0(N) \to A_\psi \to E. \]
This concludes the proof.

\end{proof}

\section{Beyond arithmetic-type hyperbolic uniformisation}\label{sec:khare}

The equivalence of the three versions of Modularity shown up to now may be summarised in the following.

\begin{theorem}[Mazur, 1991]
Let $E$ be an elliptic curve over $\Q$. Then $E$ admits a dominant morphism from $X_0(N)$ defined over $\Q$ for some $N \geq 1$ if and only if $E$ admits a hyperbolic uniformisation of arithmetic type over $\C$.
\end{theorem}

Note that, as shown in Step (3) of the proof in \Cref{sec:compare_mazur}, one may replace the curve $X(N)$ in the definition of hyperbolic uniformisation of arithmetic type with the modular curve $X_1(N)$ (possibly for a different $N$).

The above framing led Serre to ask for a characterisation of elliptic curves over $\C$ which admit a hyperbolic uniformisation of arithmetic type. That is, \emph{can elliptic curves over more general subfields of $\C$ be modular in Mazur's sense?}

All such curves (by Belyi) are defined over $\overline{\Q}$; thus one wishes to characterise those elliptic curves over $\overline{\Q}$ which are quotients of some modular curve $X_1(M)$. Ribet gave such a characterisation \cite{Ribet2004}, at the time conjectural on Serre's conjecture on modularity of mod-$p$ Galois representations (now a theorem of Khare--Wintenberger \cite{kharewintenberger2009I,kharewintenberger2009II}): they are precisely the ``$\Q$-curves'', those curves isogenous to all of their Galois conjugates.

So what about elliptic curves over $\overline{\Q}$ that are not $\Q$-curves? In what sense could they be modular? If one stays in the realm of uniformisations (that is, being dominated by modular curves), one is then forced to consider just plain-old hyperbolic uniformisations, not necessarily of arithmetic type. That is, one considers modular curves $X(\Gamma)$ associated to general subgroups $\Gamma$ of finite index in $\SL_2(\Z)$.

Following Khare \cite{khare2002belyi}, we make the following definition.

\begin{definition}
    Let $E$ be an elliptic curve defined over $\overline{\Q}$.
    \begin{enumerate}
        \item A \textbf{Belyi parametrisation of $E$} is a nonconstant map $X(\Gamma) \to E$ defined over $\overline{\Q}$ which we normalise by requiring that the $0$-cusp is sent to the origin. We say that such a $\Gamma$ \textbf{covers $E$}.
        \item Given a subgroup $\Gamma$ of $\SL_2(\Z)$ of finite index its \textbf{congruence hull $\Gamma^c$} is the smallest congruence subgroup which contains $\Gamma$. The \textbf{congruence defect $cd_\Gamma$ of $\Gamma$} is the index $[\Gamma^c : \Gamma]$. 
        \item The \textbf{congruence defect $cd_E$ of $E$} is the smallest congruence defect of a finite index subgroup of $\SL_2(\Z)$ that covers $E$.
    \end{enumerate}
\end{definition}

By Belyi's theorem, such parametrisations always exist and thus $cd_E$ is a well-defined invariant associated to the $\overline{\Q}$-isomorphism class of $E$.

Much is not known about such Belyi parametrisations; and here we mention two questions to be found in \cite{khare2002belyi}.

\begin{question}
    Given a number field $K$ is there a constant $c_K$ such that for all elliptic curves defined over $K$, we have
    $cd_E \leq c_K$?
\end{question}

This is true for $K = \Q$, where $c_K = 1$. But this is not known for any other number field.

We also have from above that the class of elliptic curves with congruence defect $1$ are precisely the $\Q$-curves. We thus make the following definition.

\begin{definition}
The \textbf{$\Q$-field} of an elliptic curve over $\overline{\Q}$ is the fixed field of $\Gal(\overline{\Q}/\Q)$ which fixes the $\overline{\Q}$-isogeny class of $E$.
\end{definition}

Thus the $\Q$-field of $E$ is precisely $\Q$ when $E$ is a $\Q$-curve. A refinement of the question above is then the following.

\begin{question}
    Given a number field $K$ is there a constant $c_K'$ such that for any elliptic curve with $\Q$-field $K$, we have
    $cd_E \leq c_K'$?
\end{question}

\bibliographystyle{amsalpha}

\end{document}